\documentclass[12pt]{amsart}
\usepackage{graphicx}
\numberwithin{equation}{section} 

\newtheorem{thm}{Theorem}[section]
\newtheorem{cor}[thm]{Corollary}

\newtheorem{prop}[thm]{Proposition}

\theoremstyle{definition}
\newtheorem{defn}[thm]{Definition}
\theoremstyle{remark}

\newcommand{\bea}{\begin{eqnarray}}
\newcommand{\eea}{\end{eqnarray}}
\newcommand{\ba}{\begin{array}}
\newcommand{\ea}{\end{array}}
\newcommand{\bc}{\begin{center}}
\newcommand{\ec}{\end{center}}
\newcommand{\be}{\begin{equation}}
\newcommand{\ee}{\end{equation}}

\def\l{\lambda}

\def\a{\alpha}

\def\b{\beta}

\begin{document}

\title[Measurable bundles]
{Measurable bundles of $C^*$-dynamical systems and its applications}

\author{Inomjon Ganiev}
\address{Inomjon Ganiev\\
 Department of Mathematics\\
Tashkent Railway Engineering Institute,\\
100167, Adilhodjaev str. 1, Tashkent, Uzbekistan} \email{{\tt
ganiev1@rambler.ru}}

\author{Farrukh Mukhamedov}
\address{Farrukh Mukhamedov\\
 Department of Computational \& Theoretical Sciences\\
Faculty of Science, International Islamic University Malaysia\\
P.O. Box, 141, 25710, Kuantan\\
Pahang, Malaysia} \email{{\tt far75m@yandex.ru} {\tt
farrukh\_m@iium.edu.my}}

\begin{abstract}
In the present paper we investigate  $L_0$-valued states and Markov
operators on $ C^*$-algebras over  $L_0$. In particular, we give
representations for $L_0$-valued state and Markov operators on $
C^*$ algebras over  $L_0$, respectively,  as measurable bundles of
states and Markov operators. Moreover, we apply the obtained
representations to study certain ergodic properties of $
C^*$-dynamical systems over $L_0$. \vskip 0.3cm \noindent

\\

\end{abstract}

\maketitle

\section{Introduction}

It is known that the theory of Banach bundles stemming from the
paper \cite{vN}, where it was proved such a theory has vast
applications in analysis. For other applications of the continuous
and measurable Banach bundled we refer the reader to
\cite{Ga,Gi,Gu,K,K3}. In the present paper we are going to deal with
$C^*$-algebras over $L_0$, such algebras are considered over modulus
of measurable functions $L_0$. Such $C^*$-algebras over $L_0$ allow
to construct meaningful examples of Banach-Kantorovich spaces. Note
that the theory of such spaces already well developed (see for
example, \cite{K}).  It is known that a structural theory of
$C^*$-modules started with papers of Kaplanski \cite{Kap} who used
such objects for an algebraic approach to the theory of
$W^*$-algebras. Consideration of $C^*$-algebras, $AW^*$-algebras and
$W^*$-algebras as modulus over their centers allowed to use methods
of Boolean-valued analysis to describe several properties of the
mentioned algebras (see \cite{K2,T}).

In \cite{GC1} it was introduced the notion  $ C^*$ algebra over ring
of all measurable functions $L_0.$ In that paper, by means of the
methods of a general theory of Banach measurable bundles (see
\cite{Gu}), it has been shown that any $C^*$-algebras over $L_0$ can
be represented as a measurable bundle of $C^*$-algebras.
 One of the main results in the theory of
$C^*$-algebras is the Gelfand-Naimark's theorem, which describes
commutative $C^*$-algebras over the complex filed $\mathbb{C}$ as an
algebra of complex valued continuous functions defined on the set of
pure states of the algebra. By means of the results of \cite{GC1} in
\cite{CGK1} it has been proved an analog of the Gelfand-Naimark's
theorem for commutative $C^*$-algebras over $L_0$.  Further, in
\cite{CGK2} GNS- representation was obtained for such
$C^*$-algebras.

In the present paper we investigate  $L_0$-valued states and Markov
operators on $ C^*$-algebras over  $L_0$. In particular in section 2
and 3 we give representations for $L_0$-valued state and Markov
operators on $ C^*$ algebras over  $L_0$, respectively,  as a
measurable bundle of states and Markov operators. In section 4, we
apply the obtained representations to study certain ergodic
properties of $ C^*$-dynamical systems over $L_0$. Note that there
are many variants of ergodic theorems for $ C^*$-dynamical systems
(see for example, \cite{DS},\cite{LL}-\cite{MT},\cite{NSZ}).

\section{Preliminaries}

  Let $(\Omega,\Sigma,\lambda)$ be a measurable space with a finite complete measure
  $\lambda$  and $L_0=L_0(\Omega)$  be an algebra of equivalence classes of complex measurable
  functions on $\Omega$.

  A complex linear space $E$ is said to be normed by $L_0$ if
  there is a map  $\|\cdot\|:E\longrightarrow L_0$ such that the
  following conditions are satisfied:
  $\|x\|\geq 0; \|x\|=0\Longleftrightarrow x=0$;
  $\|\lambda x\|=|\lambda|\|x\|$; $\|x+y\|\leq\|x\|+\|y\|$, for
  any $x,y\in E,\lambda\in \mathbb{C}.$

  A pair $(E,\|\cdot\|)$ is said to be a {\it lattice normed }space (LNS)
over $L_0$. An LNS is said to be \emph{d-decomposable} if for every
$x\in E$ and the decomposition $\|x\|=e_{1}+e_{2}$ with
 disjont positive elements $e_{1}$ and $e_{2}$ in $L_0$ there exist
$x_1, x_2\in E$, such that $x=x_1+x_2$ with $\|x_1\|=e_{1},
\|x_2\|=e_{2}$.

Suppose that $(E,\|\cdot\|)$ is  an LNS over $L_0$. A net
$\{x_{\alpha}\}$ of elements of $E$ is said to be $(bo)$-converging
to $x\in E$ (in this case we write $x=(bo)-\lim x_{\alpha}$), if the
net $\{\|x_{\alpha}-x\|\}$ $(o)$-converging to zero in $L_0$ (in
this case we write $(o)-\lim\|x_{\alpha}-x\|=0$). A net
$(x_{\alpha})_{\alpha\in A}$ is called $(bo)$- \emph{fundamental} if
$(x_{\alpha}-x_{\beta})_{(\alpha,\beta)\in A\times A}$
$(bo)$-converges to zero.

A LNS is called $(bo)$- \emph{complete} if  every $(bo)$-
fundamental net $(bo)$-converges in it. A \emph{Banach-Kantorovich
space (BKS) over $L_0$} is a $(bo)$-complete d-decomposable LNS over
$L_0$. It is well known \cite{K} that every BKS over $L_0$ admits an
$L_0$- module structure such that $\|\lambda u\|=|\lambda|\|u\|$ for
every $\lambda \in L_0,\, u\in E,$ where $|\lambda|$ is the module
of a function $\lambda \in L_0$ (see \cite{Gu,K}).

Let $\mathcal{U}$ be any  complex algebra and $\mathcal{U}$ is a
module over $L_0$ such that $(\lambda u)v=\lambda(uv)=u(\lambda v)$
for all $\lambda\in L_0,\, u,v\in\mathcal{U}$. We will consider in
$\mathcal{U}$  some $L_0$-vaued norm $\|\cdot\|$, which admits an
$L_0$- module structure BKS in $\mathcal{U}$, in particully
$\|\lambda u\|=|\lambda|\|u\|$ for all $\lambda \in L_0,\,
u\in\mathcal{U}$.

\begin{defn} \cite{CGK1,CGK2} The algebra  $\mathcal{U}$ is
called a \emph{Banach-Kantorovich
  algebra} over $L_0$ if $\|u\cdot v\|\leq\|u\|\|v\|$ for all
$u,v\in\mathcal{U}.$
\end{defn}

If $\mathcal{U}$ is  a Banach-Kantorovich
  algebra over $L_0$ with the unit $e$ such that $\|e\|=\mathbf{1},$
  then $\mathcal{U}$ is called an {\it unital Banach-
  Kantorovich}
  algebra over $L_0$, where $\mathbf{1}$ is unit element $L_0$.

  Let $\mathcal{U}$ be an $*$-algebra and Banach-Kantorovich
  algebra over $L_0$, such that $(\lambda u)^{\ast}=\bar{\lambda}u^{\ast}$
  for all
$\lambda\in L_0,\, u\in\mathcal{U}.$

\begin{defn}\cite{GC1} An algebra $\mathcal{U}$ is called
$C^{*}$\emph{-algebra over}
  $L_0$ if  for every  $u\in\mathcal{U}$  the equality
        $\|u\|^{2}=\|u^{*}\cdot u\|$ holds.
\end{defn}

{\sc Example} \cite{CGK2}. Let us provide an example of
$C^*$-algebra over $L_0$. To do it, let us recall some definitions
taken from \cite{CGK2}. Consider a modulus $\mathcal{A}$ over $L_0$.
A mapping
$\langle\cdot,\cdot\rangle:\mathcal{A}\times\mathcal{A}\mapsto L_0$
is called {\it $L_0$-valued inner product}, if for every
$x,y,z\in\mathcal{A}$, $\l\in L_0$ one has $\langle x,x\rangle\geq
0$; $\langle x,x\rangle=0$ if and only if $x=0$; $\langle
x,y\rangle=\overline{\langle y,x\rangle}$;  $\langle \l
x,y\rangle=\l\langle y,x\rangle$; $\langle x+y,z\rangle=\langle
x,z\rangle+\langle y,z\rangle $.

If $\langle\cdot,\cdot\rangle:\mathcal{A}\times\mathcal{A}\mapsto
L_0$ is a {\it $L_0$-valued inner product}, the formula
$\|x\|=\sqrt{\langle x,x\rangle}$ defines a $d$-decomposable
$L_0$-valued norm on $\mathcal{A}$. Then the pair
$(\mathcal{A},\langle\cdot,\cdot\rangle)$ is called {\it
Hilbert-Kaplansky modules}, if $(\mathcal{A},\|\cdot\|)$ is BKS over
$L_0$. By the same way one can define Hilbert-Kaplansky modules over
$L^\infty$ \cite{K}.

Let $E$ and $F$ be BKS over $L_0$ (resp. $L^\infty$). An operator
$T:E\to F$ is called {\it $L_0$ (resp. $L^\infty$)-linear}, if one
has $T(\a x+\b y)=\a T(x)+\b T(y)$ for every $x,y\in E$, $\a,\b\in
L_0$ (resp. $\a,\b\in L^\infty$). A linear operator $T$ is called
{\it $L_0$-bounded} (resp. {\it $L^\infty$-bounded}) if there exists
$c\in L_0$ (resp. $c\in L^\infty$) such that $\|T(x)\|\leq c\|x\|$
for every $x\in E$. For $L_0$-linear and $L_0$-bounded operator $T$
one defines $\|T\|=\sup\{\|T(x)\|:\ \|x\|\leq \mathbf{1}\}$, which
is a norm of $T$ (see \cite{K}).

Now let $\mathcal{A}$ be Hilbert-Kaplansky modulus over $L_0$. Then
$\mathcal{A}_b=\{x\in\mathcal{A}: \ \|x\|\in L^\infty\}$ is a
Hilbert-Kaplansky modulus over $L^\infty$. By $B(\mathcal{A})$
(resp. $B(\mathcal{A}_b)$ we denote the set of $L_0$-linear,
$L_0$-bounded (resp. $L^\infty$-linear,$L^\infty$-bounded) operators
on the Hilbert-Kaplansky modules $\mathcal{A}$ over $L_0$ (resp.
$\mathcal{A}_b$). For each operator $T\in B(\mathcal{A}_b)$ there is
a conjugate operator $T^*\in B(\mathcal{A}_b)$ satisfying the
following equality
\begin{equation}\label{conj}
\langle T(x),y\rangle=\langle x,T^*(y)\rangle, \ \ \
x,y\in\mathcal{A}_b.
\end{equation}
Moreover, one has $\|T^*\|=\|T\|$, $\|T^*T\|=\|T\|^2$, and therefore
$B(\mathcal{A}_b)$ is a $C^*$-algebra over $L^\infty$ (see
\cite{K}).

One can show that for each $T\in B(\mathcal{A})$ there exists
$T^*$ in $B(\mathcal{A})$ which satisfies \eqref{conj} for all
$x,y\in \mathcal{A}$ (see \cite{CGK2}). Since each operator $T\in
B(\mathcal{A})$ can be represented as $(bo)$-limit of the sequence
of operators $T_n\in B(\mathcal{A})$ with $\|T\|\in L^\infty$,
then $B(\mathcal{A})$ is a $(bo)$-completion of
$B(\mathcal{A}_b)$. This
means that $B(\mathcal{A})$ is a $C^*$-algebra over $L_0$. \\

        Let $X$ be a mapping, which sends every point $\omega \in
        \Omega$ to some $C^{*}$-algebra
        $(X(\omega),\|\cdot\|_{X(\omega)})$. In what follows, we
        assume that $X(\omega)\neq\{0\}$ for all $\omega\in\Omega$.

        A function $u$ is said to be a \emph{section} of $X$, if it is
        defined almost everywhere in $\Omega$ and takes its value $u(\omega)\in
        X(\omega)$ for $\omega\in dom (u),$ where $dom (u)$ is the
        domain of $u$.

      Let $L$ be some set of sections.

\begin{defn} \cite{GC1} A pair $(X ,L)$ is said to be a
\emph{measurable bundle of $C^{*}$-algebras} over $\Omega,$ if
\begin{enumerate}
\item[1.] $\lambda_1 c_1+\lambda_2 c_2 \in L$ for all
$\lambda_1,\lambda_2\in \mathbb{C}$
  and $c_1,c_2\in L$,  where $\lambda_1 c_1+\lambda_2 c_2 :\omega\in
   dom (c_1)\cap  dom (c_2)\rightarrow\lambda_1 c_1(\omega)+\lambda_2
   c_2(\omega)$;

\item[2.] the function $\|c\|: \omega \in dom (c) \rightarrow\|c(\omega)\|
 _{X(\omega)}$
is measurable for all $c\in L$;

\item[3.] for every $\omega\in\Omega,$ the set $\{c(\omega):c\in L,
 \omega\in  dom (c)\}$
 is dense in $X(\omega)$.

\item[4.] if $c\in L$ then $c^{*}\in L$ where $c^{*}:\omega\in dom
(c)\rightarrow c(\omega)^{*};$

\item[5.]  if $c_1,c_2\in L$, then $c_1\cdot c_2\in L$ where  $c_1\cdot
c_2:\omega\in dom (c_1)\cap  dom (c_2)\rightarrow
   c_1(\omega)\cdot c_2(\omega)$.
   \end{enumerate}
\end{defn}

   A section $s$ is said to be \emph{simple,}
    if there are $c_i\in L,A_i\in \Sigma,i=\overline{1,n}$ such that
     $s(\omega)=\sum\limits_{i=1}^n \chi_{A_i}(\omega)c_i(\omega)$
     for almost all $\omega\in\Omega.$

    A section  $u$ is called \emph{measurable} if there is a sequance
    $\{s_n\}$ of simple sections such that $s_n(\omega)\rightarrow u(\omega)$
   for almost all  $\omega\in\Omega$.

   The set of all measurable sections is denoted by $M(\Omega,X),$
   and $L_0(\Omega,X)$ denotes the factorization of this set with
   respect to equality everywhere.  We denote by $\widehat{u}$ the
   class from $L_0(\Omega,X)$ containing a section $u\in
   M(\Omega,X)$, and by $\|\widehat{u}\|$ the element of $L_0$
   containing the function $\|u\|_{X(\omega)}.$

For $\hat{u}\in L_0(\Omega,X)$ we put
$\hat{u}^{\ast}=\widehat{u(\omega)^{\ast}}.$

\begin{thm}\label{1.4} \cite{GC1} Let the pair $(X ,L)$ be a measurable bundle
of $C^{*}$-algebras, then $L_0(\Omega,X)$ is a $C^{*}$-algebra over
$L_0.$
\end{thm}

The set of all bounded measurable functions on $\Omega$ will be
denoted by $\mathcal L^{\infty}(\Omega),$ with the norm
$\|f\|_{\mathcal L^{\infty}(\Omega)}=\sup\limits_{\omega\in
\Omega}|f(\omega)|.$ Let $L^{\infty}(\Omega)=\{\widehat{f}\in L_0:
\exists \alpha>0, |\widehat{f}|\leq \alpha\mathbf{1}\}$ with the
norm $\|\widehat{f}\|_{L^{\infty}(\Omega)}= \inf\{\alpha>0:
|\widehat{f}|\leq \alpha\mathbf{1}\}.$

Let ${\mathcal L^{\infty}}(\Omega,X)=\{u\in
M(\Omega,X):\|u(\omega)\|_{X(\omega)}\in \mathcal
L^{\infty}(\Omega)\}$ and $L^{\infty}(\Omega,X)=\{\widehat{u}\in
L_0(\Omega,X): \|\widehat{u}\|\in L^{\infty}(\Omega)\}.$  One can
define the spaces $\mathcal {L^{\infty}}(\Omega,X)$ and
$L^{\infty}(\Omega,X)$ with real-valued norms $\|u\|_{\mathcal
L^{\infty}(\Omega,X)}=\sup\limits_{\omega\in
\Omega}|u(\omega)|_{X(\omega)}$ and
$\|\widehat{u}\|_{\infty}=\bigg\|\|\widehat{u}\|\bigg\|_{L^{\infty}(\Omega)},$
respectively.

It is known \cite{K},\cite{Gu} that there is a homomorphism
 $p:L^{\infty}(\Omega)\rightarrow\mathcal{L^{\infty}}(\Omega)$
 being
 a lifting  such that
\begin{enumerate}
\item[1.] $p(\widehat{f})\in \widehat{f}$ and $dom
p(\widehat{f})=\Omega;$

\item[2.] $\|p(\widehat{f})\|_{\mathcal
L^{\infty}(\Omega)}=\|\widehat{f}\|_{L^{\infty}(\Omega)}.$
\end{enumerate}

 The homomorphism $p$ is usually called a \emph{lifting} from
$L^{\infty}(\Omega)$ to $\mathcal{L^{\infty}}(\Omega).$

\begin{defn}\cite{GC1} The map
$\ell:L^{\infty}(\Omega,X)\rightarrow
     {\mathcal{L}}^{\infty}(\Omega,X)$
     is called \emph{a vector-valued lifting }(associated with $p$), if
     for all
      $\hat{u},\hat{v}\in L^{\infty}(\Omega,X)$ and $\lambda\in
      L^{\infty}(\Omega)$ the following conditions are valid:
\begin{enumerate}
\item[1.]  $\ell(\hat{u})\in \hat{u},\, dom\,
      \ell(\hat{u})=\Omega$;

\item[2.]  $\|\ell(\hat{u})(\omega)\|_{X(\omega)}=p(\|\hat{u}\|)(\omega)$;

     \item[3.]  $\ell(\hat{u}+\hat{v})=\ell(\hat{u})+\ell(\hat{v})$;

     \item[4.]  $\ell(\lambda \hat{u})=p(\lambda)\ell(\hat{u})$;

     \item[5.]  $\ell(\hat{u}^{*})=\ell(\hat{u})^{*}$;

    \item[6.]  $\ell(\hat{u}\hat{v})=\ell(\hat{u})\ell(\hat{v})$;

     \item[7.]  for every $\omega\in \Omega$ the set $\{ \ell(\hat{u})(\omega): \hat{u}\in L^{\infty}(\Omega,X)\}$ is dense in $
      {X(\omega)}.$
\end{enumerate}
\end{defn}

\begin{thm}\label{1.6}\cite{GC1} For any $C^{*}$-algebra $\mathcal{U}$
over $L_0,$ there exist measurable bundle of $C^{*}$-algebras $(X
,L)$ with vector-valued lifting such that $\mathcal{U}$ is
isometrically and *-isoomomorphic to $L_0(\Omega,X).$
\end{thm}

Let us consider a unital  $C^{*}$-algebra $\mathcal{U}$ over
 $L_0,$ in this case due to above given Theorem \ref{1.6}, we may identify
 $\mathcal{U}$ with $L_0(\Omega,X).$

Recall that a mapping $f:\mathcal{U}\rightarrow L_0$ is called
\emph{$L_{0}$-linear}, if  $f(\alpha x +\beta y )=\alpha f(x)+\beta
f(y)$ for all $\alpha,\beta \in L_0 ,x , y \in \mathcal{U}.$
 An $L_{0}$-linear functional $f:\mathcal{U}\rightarrow L_0$ is called \emph{$L_0$-bounded},
 if there exists $c\in L_{0}$
such that $\|f(x)\|\leq c\|x\|$ for all $x\in \mathcal{U}.$ For
$L_{0}$-linear $L_0$-bounded functional $f:\mathcal{U}\rightarrow
L_0$ we put
$\|f\|=\sup\{|f(x)|:x\in\mathcal{U},\,\|x\|\leq\textbf{1}\}$. An
$L_{0}$-linear  functional $f:\mathcal{U}\rightarrow L_0$ is said
be:
 \emph{positive} ($f\geq0$),  if $f(xx^{*})\geq0$ for all
$x\in\mathcal{U}$; \emph{a $L_{0}-$state}, if $f\geq0$ and
$\|f\|=\textbf{1}$ (see \cite{CGK2}).

Let  $\varphi$ be a positive $L_0$-linear functional, and $,a,b\in
\mathcal{U},\,\lambda\in L_0$. Then one has $\varphi((\lambda
a+b)^{*}(\lambda a+b))\geq0,$ which implies
$$|\lambda|^{2}\varphi(a^{*}a)+\overline{\lambda}\varphi(a^{*}b)+
 \lambda\varphi(b^{*}a)+\varphi(b^{*}b)\geq0.$$ Therefore
 $$
 \varphi(a^{*}b)=\overline{\varphi(b^{*}a)},\quad
|\varphi(a^{*}b)|^{2}\leq\varphi(a^{*}a)\varphi(b^{*}b).\eqno (1)
$$
Consequently, for positive $L_0$-linear functional $\varphi$ we have
we have $\varphi(a^{*})=\overline{\varphi(a)}$. Moreover, if
$\varphi(e)=0,$ then $\varphi=0.$

The following proposition contains properties of positive
functionals of $C^{\ast}$-algebra  over $L_0$ which are similar to
ones defined on $C^{\ast}$-algebras.

\begin{prop}\label{1.7} \cite{CGK1}  Let $\mathcal{U}^{\ast}$ be  the
set of all $L_0$-bounded $L_0$-linear functionals on $\mathcal{U}$.
Then
\begin{enumerate}
\item[(i)] if $\varphi\geq 0$, then
$|\varphi(x)|^{2}\leq\varphi(e)\varphi(x^{*}x)\leq\varphi(e)^{2}\|x\|^{2}$.
In particular, one has $\varphi\in \mathcal{U}^{*}$ and
$\|\varphi\|=\varphi(e)$;

\item[(ii)] if $\varphi\in \mathcal{U}^{*}$ and
$\|\varphi\|=\varphi(e)$, then $\varphi\geq0$;

\item[(iii)] if $\varphi\in \mathcal{U}^{*}$ and
$\|\varphi\|=\textbf{1}=\varphi(e)$, then $\varphi$ is a
$L_0$-state;

\item[(iv)] if $\varphi,\psi\geq0$ and $\alpha, \beta \in L_0,
\alpha, \beta\geq0$, then $\alpha\varphi+\beta\psi\geq0$
 and
$\|\alpha\varphi+\beta\psi\|=\alpha\|\varphi\|+\beta\|\psi\|$, in
particular, the set  $E_{\mathcal{U}}$ of all $L_0$-states on
$\mathcal{U}$ is a convex set.
\end{enumerate}
\end{prop}
\section{Measurable bundles of states}

Let  $(X,L)$  be a measurable bundle of $C^*$-algebras  and
$\varphi_{\omega}$ be a state on $X(\omega)$, for all
$\omega\in\Omega$, respectively and $M(\Omega)$ the set measurable
complex functions on $\Omega$.

\begin{defn}  The family $\{\varphi_{\omega}\}$ is
called a \textit{measurable bundle of states}, if
$$\varphi_{\omega}(x(\omega))\in M(\Omega)$$ for all $ x\in
M(\Omega,X).$
\end{defn}

\begin{thm}\label{2.2} Let $\{\varphi_{\omega}\}$ be a measurable
bundle of states.  Then a linear mapping
$\widehat{\varphi}:L_0(\Omega,X)\rightarrow L_0$  defined by
$$\widehat{\varphi}(\widehat{x})=\widehat{\varphi_{\omega}(x(\omega))}$$ is
a $L_0$-state on $L_0(\Omega,X)$.
\end{thm}

\begin{proof} From \cite{GK} we find that  $\widehat{\varphi}$ is an
$L_0$ -linear, $L_0$ -bounded operator on $L_0(\Omega,X)$. According
to $\varphi_{\omega}(x(\omega)x(\omega)^*)\geq 0$ for a.e.
$\omega\in\Omega$, one gets
$\widehat{\varphi}(\widehat{x}\widehat{x}^*)\geq 0$, and therefore
one has  $\|\widehat{\varphi}\|=\widehat{\varphi}(\widehat{e}).$ On
the other hand, we have
$\widehat{\varphi}(\widehat{e})=\widehat{\varphi_{\omega}(e(\omega))}=\mathbf{1}$.
Therefore
$\|\widehat{\varphi}\|=\widehat{\varphi}(\widehat{e})=\mathbf{1}.$
 Then by Proposition 2.7 $\widehat{\varphi}$ is an $L_0$-
state.
\end{proof}

\begin{thm}\label{2.3} Let $\varphi$ be a $L_0$-state on
$L_0(\Omega,X)$.
 Then there exists a measurable bundle of states $\{\varphi_{\omega}\}$ such that
 $$\varphi(x)(\omega)=\varphi_{\omega}(x(\omega))$$
 for a.e.
$\omega\in\Omega$ and for any $x\in L_0(\Omega,X)$.
\end{thm}

\begin{proof} Since $\|{\varphi}\|={\varphi}({e})=\mathbf{1},$
we have $\varphi(x)\in L^{\infty}(\Omega)$
  for any $x\in L^{\infty}(\Omega,X)$.

Define  a linear functional $\varphi_{\omega}$ on
$\{\ell(x)(\omega): x\in L^{\infty}(\Omega,X)\}$
  by
$$\varphi_{\omega}(\ell(x)(\omega))=p(\varphi(x))(\omega)$$
where $p$ is a lifting on $ L^{\infty}(\Omega).$ Due to
\begin{eqnarray*}
|\varphi_{\omega}(\ell(x)(\omega))|&=&|p(\varphi(x))|(\omega)\\
&=& p(|\varphi(x)|)(\omega)\\
&\leq&
  p(\|x\|)(\omega)
\end{eqnarray*}
we find that $\varphi_{\omega}$ is bounded and correctly defined.
   Using positivity $p$ and $\varphi,$ we have
$$\varphi_{\omega}(\ell(x)(\omega)\ell(x)^*(\omega))=
\varphi_{\omega}(\ell(xx^*)(\omega))=p(\varphi(xx^*))(\omega)\geq
0.$$ Therefore  $\varphi_{\omega}$ is a positive linear form on
$\{\ell(x)(\omega): x\in L^{\infty}(\Omega,X)\}$. Using the density
of $\{\ell(x)(\omega): x\in L^{\infty}(\Omega,X)\}$ in $X(\omega)$
for any $\omega\in\Omega,$ by means of the continuity argument we
can extend the linear form $\varphi_{\omega}$ as follows:
$$\varphi_{\omega}(x(\omega))=\lim\limits_{n\rightarrow\infty}\varphi_{\omega}(\ell(x_n)(\omega)),$$
where $x_n\in L^{\infty}(\Omega,X).$ The extension is also denoted
by $\varphi_{\omega}$ which is clearly positive.

We put $e(\omega)=\ell(e)(\omega)$ for all $\omega\in\Omega.$
Since $e$ is unit element in $L_0(\Omega,X)$ we have that
$e(\omega)$ is unit element in $X(\omega).$

We are going to show that $\varphi_{\omega}$ is a state in
$X(\omega)$. Inded, let $x(\omega)\in X(\omega).$  Then
$$\varphi_{\omega}(x(\omega)x(\omega)^*)=
 \lim\limits_{n\rightarrow\infty}\varphi_{\omega}(\ell(x_n)(\omega)\ell(x_n^*)(\omega)).$$

By Theorem  2.1.4. \cite{D} we derive
$\|\varphi_{\omega}\|=\varphi_{\omega}(e(\omega))$.

On the other hand, from
$$\varphi_{\omega}(e(\omega))=p(\varphi(e))(\omega)=p(\mathbf{1})(\omega)=1.$$
one gets
 $\|\varphi_{\omega}\|=\varphi_{\omega}(e(\omega))=1$ for
any $\omega\in\Omega$. Hence $\varphi_{\omega}$ is a state in
$X(\omega)$.

Due to $\varphi_{\omega}(\ell(x)(\omega))\in \mathcal
L^{\infty}(\Omega),$ we obtain $\varphi_{\omega}(x(\omega))\in
M(\Omega)$ for any  $x\in M(\Omega,X).$ Therefore,
$\{\varphi_{\omega}\}$ is a measurable bundle of states.

It is clear that $\varphi_{\omega}(x(\omega))=\varphi(x)(\omega)$
for $x\in L^{\infty}(\Omega,X)$ for almost all $\omega\in\Omega$.

Let $x\in L_0(\Omega,X)$. Since $L^{\infty}(\Omega,X)$ is
$(bo)$-dence in $ L_0(\Omega,X)$, so there is a sequence  $x_n\in
L^{\infty}(\Omega,X)$  such that
$\|x_n-x\|\stackrel{(o)}{\rightarrow}0.$ Then
$\|x_n(\omega)-x(\omega)\|_{X(\omega)}\rightarrow 0$ for almost all
$\omega\in\Omega$. From
$$\varphi(x)=(o)-\lim\limits_{n\rightarrow\infty}\varphi(x_n)$$ we get
$$\|\varphi_{\omega}(x_n(\omega))-\varphi(x)(\omega)\|_{X(\omega)}=
\|\varphi(x_n)(\omega)-\varphi(x)(\omega)\|_{X(\omega)}\rightarrow
0$$
 for almost all $\omega\in\Omega$. Therefore
$\varphi(x)(\omega)=\lim\limits_{n\rightarrow\infty}\varphi_{\omega}(x_n((\omega)))$
for almost all $\omega\in\Omega$. On the other hand, the continuity
of
 $\varphi_{\omega}$ yields that
$\lim\limits_{n\rightarrow\infty}\varphi_{\omega}(x_n(\omega))=\varphi_{\omega}(x(\omega))$.
Hence for every $x\in L_0(\Omega,X)$ we have
$\varphi(x)(\omega)=\varphi_{\omega}(x(\omega))$ for almost all
$\omega\in\Omega$.
\end{proof}

\section{Measurable bundles of Markov
operators}

In this section we going to define a notion of Markov operator on
$C^*$-algebra over $L_0$, and its measurable bundle.

First recall that a linear mapping $T$ from a unital $C^*$-algebra
over $L_0$ $\mathcal{U}$ with unit $e$ into itself is called a
{\it Markov operator}, it $T$ is positive (i.e. $Tx\geq 0$
whenever $x\geq 0$, $x\in \mathcal{U}$), and $Te=e$. Similarly, an
$L_0$-linear operator $T:L_0(\Omega,X)\rightarrow L_0(\Omega,X)$
is called {\it positive} if $Tx\geq 0$ whenever $x\geq 0$, $x\in
L_0(\Omega,X)$.

The next theorem is an analog of the criterion of positivity of
operators on $C^*$-algebras for $C^*$-algebras over $ L_0(\Omega).$

\begin{thm}\label{3.1} Let $T:L_0(\Omega,X)\rightarrow
L_0(\Omega,X)$ be an $L_0$-linear operator with $Te=e.$ Then $T$ is
positive if and only if the $\|T\|=\mathbf{1}$.
\end{thm}

\begin{proof} "only if" part. Assume that $\|T\|=\mathbf{1}$. We consider the set
$L^{\infty}(\Omega,X) (\subset L_0(\Omega,X))$ with the
real-valued norm $\|x\|_{\infty}=\| \|x\|
\|_{L^{\infty}(\Omega)}$. According to \cite{CGK2} the pair
$(L^{\infty}(\Omega,X),\|\cdot\|_{\infty})$ is a $C^*$- algebra.

  Let $x\in L^{\infty}(\Omega,X).$ Then
 $\|Tx\|\leq \|x\|\in L^{\infty}(\Omega)$ implies  $T(x)\in
L^{\infty}(\Omega,X).$ Hence  we have
$$\|Tx\|_{\infty}\leq
\|x\|_{\infty}, \|T\|_{\infty}\leq 1.$$  Due to $Te=e$ and
$\|Te\|_{\infty}= \|e\|_{\infty}$ one finds $\|T\|_{\infty}=1.$ It
follows from Corollary 3.2.6 \cite{BR}  that $T\geq 0$ in
$L^{\infty}(\Omega,X)$.

Now let $x\in L_0(\Omega,X)$. Then there is a sequence
$\{x_n\}\subset L^{\infty}(\Omega,X)$ such that
$\|x_n-x\|\stackrel{(o)}{\rightarrow}0.$ Then
$\|x_nx_n^*-xx^*\|\stackrel{(o)}{\rightarrow}0$. Therefore, the
inequality
$$\|T(x_nx_n^*)-T(xx^*)\|\leq \|x_nx_n^*-xx^*\|$$
implies $T(x_nx_n^*)\stackrel{(bo)}\rightarrow T(xx^*).$ As
$T(x_nx_n^*)\geq 0,$ we get $T(xx^*)\geq 0$. Hence  $T\geq0$.

"If" part. Now assume that $T\geq0$. For any $x\in L_0(\Omega,X)$ we
define an element $\alpha(\omega)$ as follows:

$$\alpha(\omega)=\left\{%
\begin{array}{ll}
    0 & \hbox{if $\|x\|(\omega)=0$;} \\[2mm]
    \frac{1}{\|x\|(\omega)} & \hbox{if $\|x\|(\omega)\neq0$.} \\
\end{array}%
\right.
$$

   Put $z=\alpha x$, then $\|z\|\leq \mathbf{1}.
$  Therefore   $\|z\|_{\infty}\leq1.$ Since $T\geq0,$ then Corollary
3.2.6. \cite{BR} implies that $\|T(z)\|_{\infty}\leq1.$ So,
$\|T(z)\|\leq\textbf{1}.$ According to the construction, $x=\|x\|z$,
and therefore $T(x)=\|x\|T(z).$ The last equality yields that
$$\|T(x)\|=\|x\|\|T(z)\|\leq\|x\|.
$$
So, $\|T\|\leq\textbf{1}.$ From $T(e)=e$ and
$\|T(e)\|=\|e\|=\textbf{1},$ we get $\|T\|=\textbf{1}.$
\end{proof}

A positive $L_{0}$-linear operator $T:L_0(\Omega,X)\rightarrow
L_0(\Omega,X)$ with $Te=e$ is called a {\it Markov operator}.

\begin{defn} A family of operators
$\{T_{\omega}:X(\omega)\rightarrow X(\omega):\omega\in\Omega:
\omega\in \Omega\}$ is called a {\it measurable bundle of operators}
if $T_{\omega}u(\omega)\in M(\Omega,X)$ for all  $u(\omega)\in
M(\Omega,X)$.
\end{defn}

\begin{defn} A measurable bundle of operators
$\{T_{\omega}:X(\omega)\rightarrow X(\omega):\omega\in\Omega\}$
called a {\it measurable bundle of Markov operators} if $T_{\omega}$
is a Markov operator for almost all $\omega\in\Omega$.
\end{defn}

\begin{thm}\label{3.5} Let $\{T_{\omega}:X(\omega)\rightarrow
X(\omega):\omega\in\Omega\}$ be a measurable bundle of Markov
operators, then $L_{0}$-linear operator
$\widehat{T}:L_{0}(\Omega,X)\rightarrow L_{0}(\Omega,X)$ defined by
$$\widehat{T}\widehat{x}=\widehat{T_{\omega}x(\omega)}$$ is a Markov
operator on $L_{0}(\Omega,X)$.
\end{thm}

Proof follows from definition.

\begin{thm}\label{3.6} Let $T:L_0(\Omega,X)\rightarrow
L_0(\Omega,X)$ be a Markov operator, then there exists a measurable
bundles of Markov operators $T_{\omega}:X(\omega)\rightarrow
X(\omega)$ such that $(Tx)(\omega)=T_{\omega} x(\omega)$ for almost
all $\omega\in\Omega$ and any $x\in L_{0} (\Omega, X).$
\end{thm}

\begin{proof} From $T\geq 0$ and $Te= e,$ according to Theorem \ref{3.1} we
have $\|T\|= \textbf{1}.$ Therefore, $\|Tx\|\leq \|x\|$ for any $x
\in L_{0} (\Omega, X).$  The last inequality implies that $Tx\in
L^{\infty} (\Omega, X)$ if $x \in L^{\infty} (\Omega, X).$

Let $\ell$ be  a vector-valued lifting on $L^{\infty} (\Omega, X)$
associated with the lifting $p$. Now define a linear operator
$T_{\omega}$ from $ \{\ell(x)(\omega):x \in L^{\infty} (\Omega,
X)\}$ to $X(\omega)$ by
$$T_{\omega}(\ell(x)(\omega))=\ell(Tx)(\omega).$$
The following relations
\begin{eqnarray*}
\|T_{\omega}(\ell(x)(\omega))\|_{X(\omega)}&=&\|\ell(Tx)(\omega)\|_{X(\omega)}\\
&=&p(\|Tx\|)(\omega)\leq p(\|x\|)(\omega)\\
&=& \|\ell(x)(\omega)\|_{X(\omega)}
\end{eqnarray*}
 imply that $T_{\omega}$ is a bounded
operator and well-defined for  all $\omega\in\Omega.$ Due to the
density of the set $\{\ell(x)(\omega):x \in L^{\infty} (\Omega,
X)\}$ in $X(\omega)$ for  all $\omega\in\Omega,$ by the continuity
argument, we can extend $T_{\omega}$ to a continuous linear operator
on $X(\omega).$ This extension is also denoted by $T_{\omega}$.

It is clear that
$$T_{\omega}(e(\omega))=T_{\omega}(\ell(e)(\omega))=
\ell(Te)(\omega)=\ell(e)(\omega)=e(\omega)$$
 for  all
$\omega\in\Omega.$

Now we will show that $T_{\omega}$ is positive. Indeed, from
$$T_{\omega}(\ell(x)(\omega)\ell(x)^*(\omega))=T_{\omega}(\ell(xx^*)(\omega))=\ell(T(xx^*))(\omega)\geq 0$$
 we get $T_{\omega}\geq 0$ on
$\{\ell(x)(\omega):x \in L^{\infty} (\Omega, X)\}$ for  all
$\omega\in\Omega.$

Now let $x(\omega)\in X(\omega).$ Then there exists a sequence
$\{x_n\} \subset L^{\infty} (\Omega, X)$ such that
$\ell(x_n)(\omega)\rightarrow x(\omega)$ in the norm of $X(\omega)$
and
\begin{eqnarray*}
&&
T_{\omega}(x(\omega))=\lim\limits_{n\rightarrow\infty}T_{\omega}(\ell(x_n)(\omega)),\\
&&T_{\omega}(x(\omega)x^*(\omega))=\lim\limits_{n\rightarrow\infty}T_{\omega}(\ell(x_n)(\omega)\ell(x_n)^*(\omega)).
\end{eqnarray*}
From $T_{\omega}(\ell(x_n)(\omega)\ell(x_n)^*(\omega))\geq 0,$ we
obtain  $T_{\omega}(x(\omega)x^*(\omega))\geq 0$. This completes the
proof.
\end{proof}

\section{Measurable bundles of $C^{\ast}$-
dynamical systems}

In this section we study measurable bundles of $C^*$-dynamical
systems over $L_0$.

Recall that any triplet $(A,\varphi,T)$, consisting of a
$C^*$-algebra $A$, a state $\varphi$ on $A$ and a Markov operator
$T:A\mapsto A$ with $\varphi\circ T=\varphi$, is called {\it a state
preserving $C^*$-dynamical system}. Similarly, we define {\it a
state preserving $C^*$-dynamical system over $L_0$} as a triplet
$(\mathcal{U},\hat\varphi,\hat T)$ consisting of a $C^*$-algebra
$\mathcal{U}$ over $L_0$, an $L_0$-state $\hat\varphi$ on
$\mathcal{U}$ and a Markov operator $\hat T$ on $\mathcal{U}$ with
$\hat\varphi\circ \hat T=\hat \varphi$,

Let  $(X,L)$ is a  measurable bundle of $C^*$-algebras.

\begin{defn} A collection of state preserving $C^*$-dynamical systems $\{(X(\omega),
\varphi_{\omega},T_{\omega}): \ \omega\in\Omega\}$ is said to be a
{\it measurable bundle of $C^*$- dynamical systems}, if
\begin{enumerate}

\item[(i)] The correspondence $\omega\rightarrow \{T_{\omega}\}$ is a  measurable bundle of
Markov operators;

\item[(ii)] the correspondence $\omega\rightarrow \{\varphi_{\omega}\}$ is a
measurable bundle of states.
\end{enumerate}
\end{defn}

\begin{thm}\label{4.3} Let $(\mathcal{U}, \varphi,T)$ be a
$C^{\ast}$-dynamical system over $L_{0}$, then there exists a
measurable bundle of $C^{\ast}$-dynamical systems $(X(\omega),
T_{\omega},\varphi_{\omega})$ such that
\begin{enumerate}
\item[(i)] $\mathcal{U}$ is isometrically and *-isomorphic to
$L_0(\Omega,X)$;

\item[(ii)] $T=\widehat{T_{\omega}}$;

\item[(iii)] $\varphi=\widehat{\varphi_{\omega}}$
\end{enumerate}
\end{thm}

\begin{proof} (i) immediately follows from Theorem \ref{1.6}, and (ii) follows from Theorem
\ref{3.6}. The last (iii) follows from Theorem \ref{2.3}.
\end{proof}

We say that $(\mathcal{U},\varphi, T)$ - $C^*$-dynamical system over
$L_{0}$ is {\it ergodic} if
\begin{equation}\label{erg}
(o)-\lim
\frac{1}{n}\sum\limits_{k=0}^{n-1}\varphi(y{T}^k({x}))=\varphi(y){\varphi}({x})
\ \ \textrm{for all} \ \ x,y\in\mathcal{U}.
\end{equation}

\begin{thm}\label{erg1}  If  measurable bundles of $C^*$- dynamical systems
$(X(\omega),\varphi_{\omega}, T_{\omega})$ are ergodic for almost
all $\omega\in\Omega$, then $(\mathcal{U}, \varphi,T)$ is ergodic.
\end{thm}

\begin{proof}
The proof  is obvious.
\end{proof}

Let $(\mathcal{U}, \varphi,T)$  be $C^{\ast}$-dynamical system over
$L_{0}$. Denote
$$
\mathcal{U}^T=\{x\in\mathcal{U}: \ Tx=x\}.
$$

\begin{prop}\label{fix}  Let $(\mathcal{U}, \varphi,T)$  be $C^{\ast}$-dynamical system over
$L_{0}$, and   $(X(\omega),\varphi_{\omega}, T_{\omega})$ be its
measurable bundles of $C^*$- dynamical systems. Then   the following
assertions hold true:
\begin{enumerate}
\item[(i)] if $\widehat{x}\in \mathcal{U}^T,$ then $x\in M(\Omega,X)$
and $x(\omega)\in X(\omega)^{T_\omega}$  for almost all
$\omega\in\Omega$;

\item[(ii)] If $x(\omega)\in X(\omega)^{T_\omega}$  for almost all
$\omega\in\Omega$ and $x\in M(\Omega,X)$, then $\widehat{x}\in
\mathcal{U}^T.$
\end{enumerate}
\end{prop}

\begin{proof} (i). Let $\widehat{x}\in \mathcal{U}^T (\subset
L_0(\Omega,X))$. Then obviously $x\in M(\Omega,X)$. Since
$T\widehat{x}=\widehat{x}$, we have that
$T_{\omega}x(\omega)=(T\widehat{x})(\omega)=x(\omega)$ for almost
all $\omega\in\Omega$. This means that $x(\omega)\in
X(\omega)^{T_\omega}$ for almost all $\omega\in\Omega$.

(ii). Let $x(\omega)\in X(\omega)^{T_\omega}$ for almost all
$\omega\in\Omega$ and $x\in M(\Omega,X)$. Since
$T_{\omega}x(\omega)=x(\omega)$ for almost all $\omega\in\Omega$ we
get
$T\widehat{x}=\widehat{T_{\omega}x(\omega)}=\widehat{x(\omega)}=\widehat{x}.$
Thus $\widehat{x}\in \mathcal{U}^T.$
\end{proof}

\begin{cor}\label{erg2}
Let a measurable bundles of $C^*$- dynamical systems
$(X(\omega),\varphi_{\omega}, T_{\omega})$ with $\varphi_\omega$ is
faithful state for almost all $\omega\in\Omega$, be ergodic. Then
$(L_0(\Omega,X),\hat\varphi,\widehat T)$ is an ergodic
$C^*$-dynamical system over $L_0$. Moreover, one has
$L_0(\Omega,X)^{\widehat T}=L_0\hat e$.
\end{cor}

\begin{proof} From Theorem \ref{erg1} we immediately find that $(L_0(\Omega,X),\hat\varphi,\widehat T)$ is
ergodic. Now ergodicity of $(X(\omega),\varphi_{\omega},
T_{\omega})$ with the faithfulness of the state $\varphi_\omega$
implies that $X(\omega)^{T(\omega)}=\mathbb{C}e({\omega})$ for
almost all $\omega\in\Omega$. Hence, from Proposition \ref{fix} we
find that  $L_0(\Omega,X)^{\widehat T}=L_0\hat e$
\end{proof}

\begin{defn} Let  $(\mathcal{U},\varphi,T)$  be $C^*$-
dynamical system over $L_{0}$. Then it is called {\it uniquely
ergodic} if
\begin{equation}\label{UE1} (bo)-\lim
\frac{1}{n}\sum\limits_{k=0}^{n-1}{T}^k({x})={\varphi}({x}){e}.
\end{equation}
\end{defn}

From this definition we immediately find that unique ergodicity of
$(\mathcal{U},\varphi,T)$ -$C^*$- dynamical system over $L_{0}$
implies its ergodicity. Moreover, in this case, we have
$\mathcal{U}^T=L_0e$.

\begin{thm}\label{UE-1}  Let  $(L_0(\Omega,X), \hat{\varphi},\widehat{T})$  be $C^{\ast}$-dynamical system over
$L_{0}$, and   $(X(\omega),\varphi_{\omega}, T_{\omega})$ be its
measurable bundles of $C^*$- dynamical systems. If
$(X(\omega),\varphi_{\omega}, T_{\omega})$ is uniquely ergodic for
almost all $\omega\in\Omega$, then $(L_0(\Omega,X),
\hat{\varphi},\widehat{T})$ is uniquely ergodic.
\end{thm}

\begin{proof} Let measurable bundle of $C^*$-- dynamical systems
$(X(\omega),\varphi_{\omega}, T_{\omega})$ be uniquely ergodic for
almost all $\omega\in\Omega.$ According to Theorem 3.2 \cite{MT} we
get
$$\frac{1}{n}\sum\limits_{k=0}^{n-1}T_{\omega}^k(x(\omega))\rightarrow
 \varphi_{\omega}(x(\omega))e(\omega)$$
in the norm of $X(\omega)$ for almost all $\omega\in\Omega.$

From the equality
$$\widehat{\bigg\|
\frac{1}{n}\sum\limits_{k=0}^{n-1}T_{\omega}^k(x(\omega))-
 \varphi_{\omega}(x(\omega))e(\omega) \bigg\|_{X(\omega)}}=
 \bigg\|\frac{1}{n}\sum\limits_{k=0}^{n-1}\widehat{T}^k(\widehat{x})-
 \varphi(\widehat{x})\widehat{e}\bigg\|_{L_{0}(\Omega, X)}$$
  we have
 $$\frac{1}{n}\sum\limits_{k=0}^{n-1}\widehat{T}^k(\widehat{x}) \stackrel{(bo)}{\longrightarrow } \varphi(\widehat{x})\widehat{e}$$
in $L_{0}(\Omega,X).$
\end{proof}

{\sc Example 1.} Let $\mathcal{U}=M_2(L_{0})=\left\{\left(
                                                           \begin{array}{cc}
                                                             \alpha_{11} & \alpha_{12} \\
                                                             \alpha_{21} & \alpha_{22} \\
                                                           \end{array}
                                                         \right)
 :\alpha_{ij}\in L_{0}\right\}$ be a martix $C^{\ast}$-algebra over
 $L_{0}$. Then  $X(\omega)=M_2(\mathbb{C}).$ Let $\mathcal{{E}}:M_2(\mathbb{C})
 \otimes M_2(\mathbb{C})\rightarrow
 M_2(\mathbb{C})$ be the canonical conditional expectation, i.e.
 $\mathcal{{E}}(x(\omega)\bigotimes y(\omega))=
 \varphi_\omega(y(\omega))x(\omega)$, where $\varphi_\omega$ is the
 canonical trace on $M_2(\mathbb{C})$.
  Take $V\in M_2(\mathbb{C})\bigotimes M_2(\mathbb{C})$ such that
  $\mathcal{E}(VV^*)=e(\omega).$ Define
  $(T_V)_\omega:X(\omega)\rightarrow X(\omega)$ by
  $$(T_V)_\omega x(\omega)=\mathcal{{E}}(V(e(\omega)\otimes
  x(\omega))V^*), \ \ \  x(\omega)\in X(\omega).
  $$ Let
  $$H=\left(
                                                  \begin{array}{cccc}
                                                    0 & 0 & 0 & 0 \\
                                                    0 & 0 & 1 & 0 \\
                                                    0 & 1 & 0 & 0 \\
                                                    0 & 0 & 0 & 0 \\
                                                  \end{array}
                                                \right)
  $$
   and $V=\sqrt{\frac{2}{1+\cosh(2\beta)}}\cdot e^{\beta H}$, then
   one can see that
  $(T^n_V)_\omega\rightarrow \varphi_\omega e(\omega)$, i.e.
  $(X(\omega),\varphi_{\omega}, (T_V)_\omega)$ is uniquely ergodic (see \cite{M1}).
  Therefore, by Theorem \ref{UE-1} one finds that $(M_2(L_{0}),\widehat{\varphi_\omega},
  T_V)$ is uniquely ergodic.

{\sc Example 2.} Let $\mathcal{U}=L_0(\Omega,C_{[0,1)})$ be the set
all measurable by Bohner vector functions with values in
$C_{[0,1)}$.
  Then with respect to $L_{0}$- valued norm
  $\|\widehat{x}\|=\|\widehat{x(\omega)}\|_{C_{[0,1)}}$ $\mathcal{U}$ is BKS
  \cite{K}, and moreover, it is
  a $C^{\ast}$-algebra  over $L_{0}$ (see \cite{CGK1}). In this case, by
  construction $X(\omega)=C_{[0,1)}$ for almost
all $\omega\in\Omega$.
  Consider on $X(\omega)$ a state $\varphi_\omega$ defined by
$\varphi_\omega(x(\omega))=\int\limits_0^1
  x(\omega,t)dt$,  and Markov operator  defined by  $T_\omega
  x(\omega,t)=x(\omega,t+\alpha(mod1))$, where $\alpha$ an irrational
  number. One can see that  $\omega\rightarrow \varphi_\omega(x(\omega))$ and
  $\omega\rightarrow T_\omega
  x(\omega,t)$ are measurable for every $x\in L_0(\Omega,C_{[0,1)}).$
  By \cite{M2,W} the $C^*$- dynamical system
  $(X(\omega),T_\omega,\varphi_\omega)$ is uniquely ergodic.

 Now let us define an $L_{0}$--state and Markov
  operator on $L_0(\Omega,C_{[0,1)})$, respectively,  by
  $$\widehat{\varphi}(\widehat{x})=\widehat{\varphi_\omega(x(\omega))}, \ \ \
  \widehat{T}(\widehat{x})=\widehat{T_\omega
  x(\omega,\cdot)}.$$
  Then by Theorem \ref{UE-1} we obtain that $(L_0(\Omega,C_{[0,1)}),\widehat{\varphi},
  \widehat{T})$ is uniquely ergodic.\\

 Let  $T:(L^{\infty}(\Omega, X),\|\cdot\|_{\infty})\rightarrow (L^{\infty}(\Omega,
 X),\|\cdot\|_{\infty})$ be a Markov operator and
 $\varphi$ be a state in $L^{\infty}(\Omega,X).$

\begin{defn} The triple
$(L^{\infty}(\Omega,X),T,\varphi)$
 is called  {\it $L^{\infty}$- uniquly ergodic} with respect to
 $\varphi$ if
 $$\frac{1}{n}\sum\limits_{k=0}^{n-1}T^k(x)\rightarrow \varphi(x)e$$
 where convergence is by the norm $\|\cdot\|_{\infty}.$
\end{defn}

 Let $T_{\omega}:X(\omega)\rightarrow X(\omega)$ be a collection
 of Markov operators and $\varphi_{\omega}$ be a collection states in
 $X(\omega).$

\begin{defn} A collection $\{T_{\omega}:
 \omega\in\Omega\}$ is
 called  {\it uniformly uniquely ergodic} with respect to
 $\varphi_{\omega}$ if
 $$\lim\limits_{n\rightarrow\infty}\sup\limits_{\omega\in\Omega}\bigg\|\frac{1}{n}\sum\limits_{k=0}^{n-1}T_{\omega}^k(x(\omega))-
 \varphi_{\omega}(x(\omega))e(\omega)\bigg\|_{X(\omega)}=0.$$
\end{defn}

 Let $$T_{\omega}(x(\omega))=\ell(Tx)(\omega)$$ and
 $$\varphi_{\omega}(x(\omega))=p(\varphi(x))(\omega)$$ for any $x \in L^{\infty}
(\Omega, X).$

\begin{thm} If $(L^{\infty}(\Omega,X),T,\varphi)$ is
$L^{\infty}$-- uniquely ergodic with respect to
 $\varphi$, then $(X(\omega),T_{\omega},\varphi_{\omega})$ is
 uniformly uniquely ergodic for all $\omega\in\Omega.$
\end{thm}

\begin{proof} If $x \in L^{\infty} (\Omega, X),$ then
\begin{eqnarray}\label{333}
\bigg\|\frac{1}{n}\sum\limits_{k=0}^{n-1}T_{\omega}^k(\ell(x)(\omega))-
 \varphi_{\omega}(\ell(x)(\omega))\ell(e)(\omega)\bigg\|_{X(\omega)}&=&
 \bigg\|\frac{1}{n}\sum\limits_{k=0}^{n-1}\ell(T^k(x))(\omega)-
 \ell(\varphi(x)e)(\omega)\bigg\|_{X(\omega)}\nonumber\\[2mm]
 &=&\bigg\|\ell\left(\frac{1}{n}\sum\limits_{k=0}^{n-1}T^k(x)-
 \varphi(x)e\right)(\omega)\bigg\|_{X(\omega)}\nonumber\\[2mm]
 &=&
 p\left(\bigg\|\frac{1}{n}\sum\limits_{k=0}^{n-1}T^k(x)-
 \varphi(x)e\bigg\|\right)(\omega)\nonumber\\[2mm]
 &\leq&
 \bigg \|p\left(\bigg\|\frac{1}{n}\sum\limits_{k=0}^{n-1}T^k(x)-
 \varphi(x)e\bigg\|\right)\bigg\|_{\mathcal{L^{\infty}}(\Omega)}\nonumber\\[2mm]
 &=&
 \bigg \|\left(\bigg\|\frac{1}{n}\sum\limits_{k=0}^{n-1}T^k(x)-
 \varphi(x)e\bigg\|\right)\bigg\|_{{L^{\infty}}(\Omega)}\nonumber\\[2mm]
 &=&
 \bigg\|\frac{1}{n}\sum\limits_{k=0}^{n-1}T^k(x)-
 \varphi(x)e\bigg\|_{{\infty}}
 \end{eqnarray}
  for all $\omega\in\Omega.$

 Therefore
 $$\sup\limits_{\omega\in\Omega}\bigg\|\frac{1}{n}\sum\limits_{k=0}^{n-1}T_{\omega}^k(x(\omega))-
 \varphi_{\omega}(x(\omega))e(\omega)\bigg\|_{X(\omega)}\leq
 \bigg\|\frac{1}{n}\sum\limits_{k=0}^{n-1}T^k(x)-
 \varphi(x)e\bigg\|_{{\infty}}.$$

 As $(L^{\infty}(\Omega,X),T,\varphi)$ is
$L^{\infty}$-- uniquely ergodic with respect to
 $\varphi$, then we get that $(X(\omega),T_{\omega},\varphi_{\omega})$ is
uniformly uniquely ergodic for all $\omega\in\Omega.$
\end{proof}

We note that the reverse of the previous theorem is also true, which
immediately follows from \eqref{333}.

\section*{Acknowledgments}
The authors are grateful to Professor Vladimir Chilin for his
valuable comments and remarks on improving the paper. The authors
acknowledge the MOHE grant FRGS11-022-0170. The second named author
(F.M.) acknowledges the Junior Associate scheme of the Abdus Salam
International Centre for Theoretical Physics, Trieste, Italy.

\end{document}